\documentstyle[12pt]{article}
\newcommand{\R}{I\!\!R}
\newcommand{\N}{I\!\!N}
\newcommand{\C}{I\!\!\!\!C}
\newcommand{\Z}{{\bf Z}}
\newtheorem{thm}{Theora}[section]
\newtheorem{Th}[thm]{Theorem}

\newtheorem{Cor}[thm]{Corollary}
\newtheorem{Lem}[thm]{Lemma}
\newtheorem{Prop}[thm]{Proposition}
\newtheorem{Def}[thm]{Definition}

\newcommand\dint{\displaystyle\int}
\newcommand\dfrac{\displaystyle\frac}
\begin{document}\title{ Real Paley-Wiener theorems \\ for the Dunkl transform on $\R^d$}
{\footnotesize
\author{Hatem MEJJAOLI   and Khalifa TRIM\`ECHE\\
{\footnotesize{\em Department of Mathematics}}\\
{\footnotesize{\em Faculty of Sciences of Tunis- CAMPUS-}}\\
{\footnotesize{\em 1060.Tunis. Tunisia.}} }}
\date{}
\maketitle \noindent{\bf Abstract}\par In this paper, we establish
real
 Paley-Wiener theorems for the Dunkl transform on $\R^d$. More
precisely, we characterize the functions in the Schwartz space
${\cal S}(\R^d)$ and in $L^2_k(\R^d)$   whose Dunkl transform has
bounded, unbounded, convex and nonconvex support.
%A characterization of
%weighted $L^2_k(U)$ spaces in terms of their images under  the
%Dunkl transformation is derived, where $U$ can be  a disc, a
%symmetric body, a nonconvex or unbounded region in $\R^d$. This
%characterization is then used to derive Paley-Wiener type theorems
%for these  spaces and characterized the support  of Dunkl
%transform of Schwartz function $f$
 %by a $L^p$ growth condition
%for all $1 \leq p \leq +\infty$.
\\

\bigskip
\section { Introduction}
\hspace*{5mm} In the last few years there has been a great
interest to real Paley-Wiener theorems  for certain integral
transforms, see \cite{Ta} for an overview references and details
for this question.\\ \hspace*{5mm} In this paper we consider the
Dunkl operators $T_j, j=1,...,d$,  which are the
differential-difference operators introduced by C.F.Dunkl in
\cite{D1}. These operators are very important in pure Mathematics
and in Physics. They  provide a useful tool in the study of
special functions with root systems (see \cite{D2}.)
\\ \hspace*{5mm}C.F.Dunkl in \cite{D3} (see also \cite{J}) has studied a Fourier
transform ${\cal F}_{D}$, called Dunkl transform defined for a
regular function $f$ by $$\forall \, x \in \R^d, \; {\cal F}_D
f(x) = \dint_{\R^d}K(-ix,y) f(y)\omega_k(y)dy, $$ where $K(-ix,y)$
represents  the Dunkl kernel and $\omega_k $ a weight function.\\
\hspace*{5mm} The aim purpose of this paper is to prove  real
Paley-Wiener theorems on
 the Schwartz space ${\cal
S}(\R^d)$ and on $L^2_k(\R^d)$. More precisely we consider first
the Paley-Wiener spaces associated with the Dunkl operators:
 $$\begin{array}{ccc}
  PW_{k}^2(\R^d) & = & \{f \in {\cal E}(\R^d)/\forall\, n\, \in\, \N,\;
   \, \triangle_k^n f \in L^2_k(\R^d) \,
   \mbox{and} \, R_f^{\triangle_k} = \displaystyle\lim_{n\to \infty}
  ||\triangle_k f||_{k,2}^{\frac{1}{2n}} < +\infty\} \\
  PW_{k}(\R^d) & = & \{f \in {\cal E}(\R^d)/\forall\, n,m\, \in\, \N, \,
  \; (1+||x||)^m \triangle_k^n f \in L^2_k(\R^d) \,
  \; \mbox{and} \, R_f^{\triangle_k}  < +\infty\},
\end{array}$$
where $ {\cal E}(\R^d) $ is the space of $C^\infty$-functions on
$\R^d$, $\triangle_k = \displaystyle\sum_{j=1}^d T_j^2$ the
Dunkl-Laplacian operator, $L^2_k(\R^d)$ the space of  square
integrable functions with respect to the measure $\omega_k(x) dx$
and $||.||_{k,2}$ the norm of the space $L^2_k(\R^d)$.\\ We
establish that ${\cal F}_D$ is a bijection from $PW_{k}^2(\R^d)$
onto $L^2_{k,c}(\R^d)$(the space of  functions in $L^2_k(\R^d)$
with compact support), and from $PW_{k}(\R^d)$ onto $D(\R^d)$(the
space of $C^\infty$-functions on $\R^d$ with compact support).\\
\hspace*{5mm}Next, we
 characterize the $L^2_k(U)$-functions by their
Dunkl transform, where $U$ is respectively  a disc, a symmetric
body, a nonconvex and an unbounded domain  in $\R^d$. These
results are the real Paley-Wiener theorems for square integrable
functions with respect to the measure $\omega_k(x)dx$.\\
\hspace*{5mm}We generalize also a theorem of H.H.Bang \cite{B} by
characterizing the support of the Dunkl transform of functions in
$ {\cal S}(\R^d) $   by an
 $L^p$ growth condition. More precisely these  real Paley-Wiener
theorems  can be stated as follow:\\ $\bullet $ The Dunkl
transform ${\cal F}_D(f)$ of $f \in {\cal S}(\R^d)$ vanishes
outside a polynomial domain $U_P = \{x \in \R^, \; P(x) \leq 1\}$,
with $P$ a non constant polynomial, if and only if $$\limsup_{n
\to +\infty}||P^n (iT)f||_{k,p} \leq 1, \; 1 \leq p \leq \infty,
$$ with $T = (T_1,...,T_d)$ and $||.||_{k,p}$ is the norm of the
space $L^p_k(\R^d)$ of $p^{th}$ integrable functions on $\R^d$
with respect to the measure $\omega_k(x)dx$.
\\ $\bullet$ A function $f \in {\cal S}(\R^d)$ is the Dunkl
transform of a function vanishing in some ball with radius $r$
centered at the origin, if and only if $$ \lim_{n \to \infty}
||\displaystyle\sum_{m = 0}^{\infty} \frac{(n\triangle_k)^m\, f
}{m!}||_{k,p}^{\frac{1}{n}}\leq \exp(-r^2), \; 1 \leq p \leq
\infty. $$ \hspace*{5mm} This paper is arranged as follows:\\
\hspace*{5mm}In the second section  we recall  the main results
about the harmonic analysis associated with the Dunkl operators.\\
\hspace*{5mm} The third section is devoted to  study the functions
such that the  support of their Dunkl transform are  compact, and
to establish the  real Paley-Wiener theorems for ${\cal F}_D$ on
 the Schawrz space ${\cal S}(\R^d)$.\\
\hspace*{5mm} In the fourth section we characterize the functions
in  ${\cal S}(\R^d)$ such that their Dunkl transform  vanishes
outside a polynomial domain.\\ \hspace*{5mm} In the fifth section
we give a necessary and sufficient  condition  for functions in
$L^2_k(\R^d)$ such that their  Dunkl transform vanishes in a disc.
\\ \hspace*{5mm} We study in the sixth section the functions such
that their Dunkl transform satisfies  the symmetric  body
property, and we derive  a real Paley-Wiener type theorem for
these functions.\\

\section{Harmonic analysis associated for the Dunkl operators.}
\hspace*{5mm} In the first two subsections we collect some
notations and results on Dunkl operators, the Dunkl kernel and the
Dunkl intertwining operators (see [6],[7],[8]).
\subsection { Reflection groups, root system and multiplicity
functions} \hspace*{5mm}We consider $\R^d$ with the euclidean
scalar product $<.,.>$ and $||x||=\sqrt{\langle x,x\rangle}$. On
${\C}^{d},\;||.||$ denotes also
the standard Hermitian norm\ while\ for\ all\ $z=(z_{1},\;...,\;z_{d}%
),\;w=(w_{1},\;...,\;w_{d})\in{\C}^{d},$%
\[
<z,w>=\displaystyle\sum_{j=1}^{d}z_{j}\overline{w}_{j}.
\]
For $\alpha\in\R^d\backslash\{0\}$, let $\sigma_{\alpha}$ be the
reflection in the hyperplan $H_{\alpha}\subset\R^d$ orthogonal to
$\alpha$, i.e.
\begin{equation}
\sigma_{\alpha}(x)=x-2\frac{\langle\alpha,x\rangle}{||\alpha||^{2}}\alpha.
\label{2.1}%
\end{equation}
A finite set $R\subset\R^d\backslash\{0\}$ is called a root system
if $R\cap\R.\alpha=\{\alpha,-\alpha\}$ and $\sigma_{\alpha}R=R$
for all $\alpha\in R$. For a given root system R the reflection
$\sigma _{\alpha},\alpha\in R$, generate a finite group $W\subset
O(d)$, the reflection group associated with R . We denote by $|W|$
its cardinality. All reflections in W correspond to suitable pairs
of roots. For a given $\beta \in\R \backslash {\alpha\in
R}{\cup}H_{\alpha}$, we fix the positive subsystem
$R_{+}=\{\alpha\in R\;/\langle\alpha,\beta\rangle>0\}$,
then for each $\alpha\in R,$ either $\alpha\in R_{+}$ or $-\alpha\in R_{+}%
$.\newline A function $k:R\longrightarrow{\C}$ on a root system
$R$ is called a multiplicity function if it is invariant under the
action of the associated reflection group $W$. If one regards $k$
as a function on the corresponding reflections, this means that k
is constant on the conjugacy classes of reflections in $W$. For
abbreviation, we introduce the index
\begin{equation}
\gamma=\gamma(k)=\displaystyle\sum_{\alpha\in R_{+}}k(\alpha). \label{2.2}%
\end{equation}
Moreover, $\omega_{k}$ denotes the weight function
\begin{equation}
\omega_{k}(x)=\prod_{\alpha\in
R_{+}}|\langle\alpha,x\rangle|^{2k(\alpha)},
\label{2.3}%
\end{equation}
which is$\;W-$invariant and homogeneous of degree
$2\gamma$.\newline We introduce the Mehta-type constant
\begin{equation}
c_{k}=(\int_{I\!\!R^{d}}\exp(-||x||^{2})\omega_{k}(x)\;dx)^{-1}, \label{2.4}%
\end{equation}

\noindent{\bf{Remark }}\\ \hspace*{5mm} For $d=1$ and
$W=\mathbf{Z}_{2}$, the multiplicity function $k$ is a single
parameter denoted $\gamma>0$ and we have $$
\forall\,x\in\R,\;\omega_{k}(x)=|x|^{2\gamma}. $$

\subsection{ Dunkl operators- The Dunkl kernel and the Dunkl intertwining
operator}

\noindent{\bf{Notations}}.  We denote by \\
 - $C(\R^{d}) (resp \;C_{c}
(\R^{d}))$\, the space of continuous functions on $\R^{d}$ (resp.
with compact support).\\
 - $C^{p}(\R^{d}) (resp \;C^{p}_{c}
(\R^{d}))$\, the space of functions of class $C^p$ on $\R^{d}$
(resp. with compact support).\\
 - $ {\cal E}(\R^{d})$ the space of
$C^{\infty}$-functions on $\R^{d}$.\\ - $ {C}^\infty_0(\R^{d})$
the space of $C^{\infty}$-functions on $\R^{d}$ which vanish at
the infinity.
\\ - ${\cal S}(\R^{d})$   the
space of $C^{\infty}$-functions  on $\R^{d}$  which are rapidly
decreasing as their derivatives.\\ - $D(\R^{d})$ the space of
$C^{\infty}$-functions  on $\R^{d}$ which are of compact
 support.\\
 We provide these spaces with the classical topology .\\\\
We consider also the following spaces\\ - ${\cal E'}(\R^{d})$ the
space of distributions on $\R^{d}$ with compact support. It is the
topological dual of ${\cal E}(\R^{d})$.\\ - ${\cal S'}(\R^{d})$
the space of tempered  distributions on $\R^{d}$. It is the
topological dual of ${\cal S}(\R^{d})$.\\

The Dunkl operators $T_{j},\; j\; = 1\;, ...,\; d $, on $\R^{d}$
associated with the
 finite reflection group W and the multiplicity function k are given by
\begin{equation}
T_{j} f(x) = \frac{\partial}{\partial x_{j}} f(x) +
\displaystyle\sum_{\alpha \in R_{+}}k(\alpha) \alpha_{j}
\frac{f(x) - f(\sigma_{\alpha}(x))}{<\alpha,x>},\quad f \; \in \;
C^{1}(\R^{d}). \label{h9}
\end{equation}
In the case $k = 0$, the $T_{j}, \, j = 1, ... , d,$ reduce to the
corresponding partial derivatives. In this paper, we will assume
throughout that $k \geq 0$ and $\gamma > 0$.\\ \hspace*{5mm} The
Dunkl Laplacian $\triangle_{k}$   on $\R^{d}$  is defined by
\begin{equation}
\triangle_{k}f = \displaystyle\sum_{j = 1}^{d}T_{j}^{2}f  =
\triangle f + 2 \displaystyle\sum_{\alpha \in R_{+}} k_{\alpha}
\delta_{\alpha} (f), \quad f \in C^{2}(\R^{d}), \label{h12}
\end{equation}
where $\triangle = \displaystyle\sum_{j = 1}^{d} \partial_{j}^{2}$
the Laplacian on $\R^{d}$ and
 $$ \delta_{\alpha}(f)(x) = \frac{<\nabla f(x),\alpha>}{<\alpha,x>} -
\frac{ f(x) - f(\sigma_{\alpha}(x))}{<\alpha,x>^{2}},$$ with
$\nabla f$  the gradient of f.\\ \hspace*{5mm}  For $f $ in   $
C_{c}^{1}(\R^{d})$ and $g \, in \, C^{1}(\R^{d})$ we have
\begin{equation}
\int_{\R^{d}} T_{j}f(x) g(x)\omega_{k}(x)\;dx = - \int_{\R^{d}}
 f(x) T_{j}g(x)\omega_{k}(x)\;dx, \, j = 1, ..., d.
\label{hh6}
\end{equation}

For $y \in \R^{d} $,  the system $$ \left\{
\begin{array}{crll}
T_{j}u(x,y) &=& y_{j} u(x,y),& j = 1, ..., d,\\\\ u(0,y) &=& 1,
&for\, all \; y \in \,\R^{d}.
\end{array}
\right. $$ admits a unique analytic solution on $\R^{d}$,  denoted
by $K(x,y)$ and called Dunkl kernel. This kernel has a unique
holomorphic extension to ${\C}^{d} \times {\C}^{d}$.\\
\noindent{\bf{Example. }}\\ \hspace*{5mm} If $d = 1$ and $W ={\
Z}_{2}$, the Dunkl kernel is given by
\begin{equation}
K(z,w) = j_{\gamma - \frac{1}{2}}(izw) + \frac{zw}{2 \gamma + 1}
 j_{\gamma + \frac{1}{2}}(izw), \quad z, \; w \in \C,
\label{h18}
\end{equation}
where for $\alpha \geq \frac{-1}{2}$, $j_{\alpha}$ is the
normalized Bessel function of index $\alpha$
 defined by
\begin{equation}
j_{\alpha}(z) = 2^{\alpha} \Gamma(\alpha + 1)
\frac{J_{\alpha}(z)}{z^{\alpha}} =  \Gamma(\alpha + 1)
\displaystyle\sum_{n = 0}^{\infty}\frac{(-1)^{n}(\frac{z}{2})^{2
n} } {n!  \Gamma(\alpha + 1 + n)}
\end{equation}
with $J_{\alpha}$ is the Bessel function of first kind and index
$\alpha$.\\\\
 \hspace*{5mm} The Dunkl kernel possesses the following
 properties

\begin{Prop}\hspace*{-2mm}.i) For all $z, w \in \C^{d}$ we have.
\begin{equation}
  K(z,w) = K(w,z) \quad ; K(z,0) = 1 \quad and\quad
K(\lambda z,w) = K(z, \lambda w),\,  for\; all \; \lambda \in \C.
\label{h20}
\end{equation}
 \hspace*{5mm} ii) For all $\nu \in \N^{d}, x \in \R^{d}$ and $z \in \C^{d}$, we have
\begin{equation}
|D_{z}^{\nu} K(x , z)| \leq ||x||^{|\nu|} \,\exp(||x|| || Rez||),
\label{h21}
\end{equation}
%\begin{equation}
%K(x , z)| \leq e^{}{||x|| || Rez||}
%\label{h22}
%\end{equation}
and for all $x, y \in \R^{d}$ :
\begin{equation}
|K(i x , y)| \leq 1, \label{h23}
\end{equation}
with $D_{z}^{\nu} = \frac{\partial^{\nu}}{\partial
z_{1}^{\nu_1}...\partial z_{d}^{\nu_d}}$ and $|\nu| = \nu_1 + ...
+ \nu_d.$\\ \hspace*{5mm} iii) For all  $x, y \in \R^{d}$ and $g
\in W$ we have
\begin{equation}
K(-i x , y) = \overline{K(i x , y)}, \quad and \quad K(g x , g y)
= K( x , y).
\end{equation}
\label{P1.2}
 \hspace*{5mm}$i\nu)$ The function $K(x,z)$ admits for all $x \in \R^{d}$ and $z \in \C^{d}$ the following Laplace type integral representation
\begin{equation}
K(x,z) = \dint_{\R^d} e^{<y,z>} d\mu_{x}(y), \label{753}
\end{equation}
where $\mu_{x}$ is a probability measure on $\R^d$, with support
in the closed ball $B(o, ||x||)$ of center o and radius
$||x||$.(See [11]).\end{Prop}

The Dunkl intertwining operator $V_k$ is defined on $C(\R^{d})$ by
\begin{equation}
\forall x \in \R^{d}, \quad V_k f(x) = \dint_{\R^d}
f(y)d\mu_{x}(y), \label{str}
\end{equation}
where $\mu_{x}$ is the measure given by the relation (\ref{753}).
\\  The operator $V_k$ satisfies the following properties
\\ \hspace*{5mm}i)We have $$ \forall x \in \R^{d}, \; \; \forall z \in \C^{d}, \;
\; K(x,z) = V_k (e^{<.,z>})(x). $$ \hspace*{5mm}  ii)The operator
$V_k$ is a topological isomorphism from ${\cal E}(\R^{d})$ onto
itself satisfying
 the transmutation relation
\begin{equation}
\forall x \in \R^{d}, \quad T_j { V}_k (f)(x) = {V}_k
(\dfrac{\partial}{\partial y_j}f)(x), \quad j = 1, ... , d,
 f \in {\cal E}(\R^d).
\label{9}
\end{equation}
\hspace*{3mm}ii) For each $x \in \R^d$ there exists a unique
distribution $\eta_x$ in ${\cal E'}(\R^d)$ with support in the
ball $B(o, ||x||)$,  such that for all f in ${\cal E}(\R^d)$ we
have
\begin{equation}
V_k^{-1}f(x) = <\eta_x, f>. \label{10}
\end{equation}
(See [16]).
\subsection{ The Dunkl transform}
\noindent{\bf{Notations}}. We denote by $L_{k}^{p}(\R^{d})$ the
space of measurable functions on $\R^{d}$ such that $$
\begin{array}{crl}
||f||_{k,p}& =& (\displaystyle \int_{\R^{d}} |f(x)|^{p}
\omega_{k}(x) \;dx)^{\frac{1}{p}} < +\infty,
 \quad if \;  1 \leq p
< + \infty,\\\\ ||f||_{k,\infty}& = & ess\; sup _{x \in \R^{d} }
|f(x)| < +\infty.
\end{array}
$$

\hspace*{5mm}The Dunkl transform of a function f in $D(\R^{d})$ is
given by
\begin{equation}
\forall y \in \R^{d}, \quad {\cal F}_{D}(f) (y) =
\dint_{\R^{d}}f(x) K(-iy,x) \omega_{k}(x)dx . \label{13}
\end{equation}
%\end{Def}
We give in the following some  properties of this transform. (See
[7][8]).
\\\\ \hspace*{5mm} i) For all f in $L_{k}^{1}(\R^{d})$ we have
\begin{equation}
||{\cal F}_{D} (f)|| _{k, \infty} \leq ||f||_{k, 1}. \label{h26}
\end{equation}
\hspace*{5mm} ii) For all $f$ in ${\cal S}(\R^{d})$ we have
\begin{equation}
\forall y  \in \R^{d}, \quad  {\cal F}_{D}( T_j f)( y ) = i y_j
{\cal F}_{D}( f ) (y) \quad , j = 1, ..., d. \label{h28}
\end{equation}
\hspace*{5mm} $iii)$ For all f in  $L_{k}^{1}(\R^{d})$ such that
${\cal F}_{D}(f)$ is in  $L_{k}^{1}(\R^{d})$, we have the
inversion formula
\begin{equation}
f(y) = \frac{c_{k}^{2}}{4^{\gamma + \frac{d}{2}}}\displaystyle
\int_{\R^{d}} {\cal F}_{D}(f)(x)
 K(i x , y) \omega_{k}(x)\; dx
, \quad a.e. \label{h31}
\end{equation}
\begin{Th}\hspace*{-2mm}. The  Dunkl transform ${\cal F}_{D}$
is a topological isomorphism. \\ \hspace*{5,5mm} i) From ${\cal
S}(\R^{d})$ onto itself. \\ \hspace*{5,5mm} ii) From $D(\R^{d})$
onto $ {H}(\C^{d})$ (the space of entire functions on $\C^{d}$,
rapidly decreasing and of exponential type.) \\\noindent The
inverse transform ${\cal F}_{D}^{- 1}$ is given by
\begin{equation}
 \forall y \in \R^{d}, \quad {\cal F}_{D}^{-1}(f)(y) = \frac{c_{k}^{2}}{4^{\gamma + \frac{d}{2}}}
{\cal F}_{D}(f)(-y), \quad f \in S(\R^{d}). \label{h32}
\end{equation}
\end{Th}
\begin{Th} \hspace*{-2mm}. i) Plancherel formula for ${\cal
F}_D$ .\\ For all f in ${\cal S}(\R^{d})$ we have
\begin{equation}
\displaystyle \int_{\R^{d}} |f(x)|^{2}
  \omega_{k}(x)\; dx = \displaystyle \frac {c_{k}^{2}} {4^{\gamma + \frac{d}{2} } }\displaystyle \int_{\R^{d}}
| {\cal F}_{D}(f)(\xi)|^{2}  \omega_{k}(\xi)\; d\xi. \label{h33}
\end{equation}
\hspace*{5mm} ii) Plancherel theorem for ${\cal F}_{D}$.\\The
renormalized Dunkl transform $f \to 2^{-( \gamma +  \frac{d}{2})}
c_{k} {\cal F}_{D}(f)$ can be uniquely extended to
 an isometric isomorphism on   $L_{k}^{2}(\R^{d})$. \label{Tp}\end{Th}
\begin{Prop}\hspace*{-2mm}.
Let $1 \leq p \leq 2$. The Dunkl transform ${\cal F}_{D}$ can be
extended to a continuous mapping from $L_{k}^{p}(\R^{d})$ into
$L_{k}^{q}(\R^{d}),$ with $q$  the conjugate component of $p$.
\label{cher}
\end{Prop}
\begin{Def}\hspace*{-2mm}. i) The Dunkl transform of a
distribution $\tau$ in ${\cal S}'(\R^{d})$ is defined by $$ <
{\cal F}_{D}(\tau), \phi > = < \tau,{\cal F}_{D }(\phi)>, \quad
\phi \in {\cal S}(\R^{d}).  $$ \hspace*{5mm} ii) We define the
Dunkl transform of a distribution $\tau$ in ${\cal E'}(\R^d)$  by
$$ \forall \, y \in \R^d, \; {\cal F}_{D}(\tau)(y) = \langle
\tau_x, K(-ix,y) \rangle.$$
\end{Def}
\begin{Th}\hspace*{-2mm}.
 The  Dunkl transform ${\cal F}_{D}$
is a topological isomorphism. \\ \hspace*{5,5mm} i) From ${\cal
S'}(\R^{d})$ onto itself. \\ \hspace*{5,5mm}ii) From ${\cal
E}'(\R^{d})$ onto ${\cal H}(\C^{d})$(the space of entire functions
on $\C^{d}$, slowly increasing and of exponential type.)
 \end{Th}
 %\begin{Prop}\hspace*{-2mm}.
\hspace*{5mm}Let $\tau$ be in ${\cal S'}(\R^d)$. We define the
distribution $T_j \tau$, $j=1,...,d,$ by $$<T_j \tau, \psi> =  -
<\tau, T_j \psi>, \; \mbox{ for\, all} \; \psi \, \in \; {\cal
S}(\R^d).$$ This distribution satisfies the following properties
\begin{eqnarray}
{\cal F}_{D}(T_j \tau) &=& i y_j {\cal F}_{D}( \tau),\quad j = 1,
..., d. \label{sol1}
\\
{\cal F}_{D}(\triangle_k \tau) &=& -|| y||^2 {\cal F}_{D}( \tau).
\label{sol2}
\end{eqnarray}
%\end{Prop}
\hspace*{5mm}We consider $f$ in $L^2_k(\R^d)$.We define the
distribution $T_f$ in ${\cal S'}(\R^d)$ by $$\langle T_f,\varphi
\rangle = \dint_{\R^d}f(x)\varphi(x)\omega_k(x)dx, \; \varphi \in
{\cal S}(\R^d).$$ In the following $T_f$  will be  denoted  by
$f$.
\begin{Prop}\hspace*{-2mm}. Let $f$ be in $L^2_k(\R^d)$. Then we have
\begin{equation}\label{ppppp}
  {\cal F}_{D}(\triangle_k f) = -||x||^2 {\cal F}_{D}(f).
\end{equation}
  \end{Prop}
\noindent{\bf{Proof}}\\\hspace*{5mm} For all $\varphi \in {\cal
S}(\R^d)$  we have $$\langle \triangle_k f,\varphi\rangle =
\langle f,\triangle_k \varphi\rangle = \dint_{\R^d}f(x)\triangle_k
\varphi(x)\omega_k(x)dx.$$ But $$\begin{array}{lll}
 \langle {\cal F}_{D}(\triangle_k f),\varphi\rangle  & = & \langle \triangle_k f,
 {\cal F}_{D}(\varphi)\rangle = \langle f,\triangle_k {\cal F}_{D}(\varphi)\rangle\\
   & = & \dint_{\R^d}f(y){\cal F}_{D}(-||x||^2 \varphi(.)) (y) \omega_k(y)dy\\
   & = &  -\dint_{\R^d}{\cal F}_{D}(f)(x)||x||^2 \varphi(x)
   \omega_k(x)dx \\ &=& \langle -|| x||^2 {\cal
   F}_{D}(f),\varphi\rangle.
\end{array}$$
Thus $${\cal F}_{D}(\triangle_k f) = -||x||^2 {\cal F}_{D}(f).$$
\noindent{\bf{Notations.}} We denote by\\
 \hspace*{5mm} - $L^2_{k,c}(\R^d)$ the space of functions in $L^2_{k}(\R^d)$
  with compact support.\\
 \hspace*{5mm} - ${\cal H}_{L^2_k}(\C^d)$ the space of entire functions $f$
 on $\C^d$  of exponential type such that $f_{|\R^d}$ belongs to  $L^2_{k}(\R^d)$.
\begin{Th}\hspace*{-2mm}.
The Dunkl transform ${\cal F}_D$ is bijective from
$L^2_{k,c}(\R^d)$ onto ${\cal H}_{L^2_k}(\C^d)$.
\end{Th}
\noindent{\bf{Proof}}\\ \hspace*{5mm}i) We consider the function
$f$ on $\C^d$ given by \begin{equation} \forall \, z \in \C^d, \;
f(z) = \dint_{\R^d}g(x) K(-ix,z)
\omega_k(x)dx,\label{159753}\end{equation} with $g \in
L^2_{k,c}(\R^d)$.\\ By derivation under the integral sign
 and by using the inequality (11), we deduce that the function
 $f$ is entire on $\C^d$ and of exponential
type.\\ On the other hand the relation (\ref{159753}) can also be
written in the form $$ \forall \; y \in \R^d, \; f(y) = {\cal
F}_{D}(g)(y).$$ Thus from Theorem \ref{Tp} the function
$f_{|\R^d}$ belongs to $L^2_{k}(\R^d)$. Thus $f \in {\cal
H}_{L^2_k}(\C^d)$.
\\ \hspace*{5mm}ii) Reciprocally let $\psi $ be in ${\cal
H}_{L^2_k}(\C^d)$. From Theorem 2.6 ii) there exists $S \in {\cal
E'}(\R^d)$ with support in the ball $B(o,a)$ of center $o$ and
radius $a$, such that
\begin{equation}\label{tgvam}
  \forall \, y \in \R^d, \; \psi(y) = \langle S_x,
  K(-ix,y)\rangle.
\end{equation}
On the other hand as $\psi_{|\R^d}$ belongs to $L^2_{k}(\R^d)$,
then from Theorem \ref{Tp} there exists \linebreak $h \in
L^2_{k}(\R^d)$ such that
\begin{equation}\label{tgvam1}
 \psi_{|\R^d} = {\cal F}_{D}(h).
\end{equation}
Thus from (\ref{tgvam}), for all $\varphi \in D(\R^d)$ we have
$$\begin{array}{lll}
  \dint_{\R^d}\psi(y) \overline{{\cal F}_{D}(\varphi)(y)}\omega_k(y)dy & = &
  \langle S_x,\dint_{\R^d}
  K(-ix,y)\overline{{\cal F}_{D}(\varphi)(y)}\omega_k(y)dy\rangle.
\end{array}$$
Thus using (22) we deduce that
\begin{equation}\label{hhhhh}
\dint_{\R^d}\psi(y) \overline{{\cal
F}_{D}(\varphi)(y)}\omega_k(y)dy  = \frac{4^{\gamma +
\frac{d}{2}}}{c_k^2}
  \langle S,\varphi\rangle.
\end{equation}
On the other hand  (\ref{tgvam1}) implies  $$\dint_{\R^d}\psi(y)
\overline{{\cal F}_{D}(\varphi)(y)}\omega_k(y)dy =
\dint_{\R^d}{\cal F}_{D}(h)(y) \overline{{\cal
F}_{D}(\varphi)(y)}\omega_k(y)dy.$$ But from Theorem 2.2 we deduce
that
\begin{equation}\label{uuuuu}\begin{array}{lll}
\dint_{\R^d}{\cal F}_{D}(h)(y) \overline{{\cal
F}_{D}(\varphi)(y)}\omega_k(y)dy &=& \frac{4^{\gamma +
\frac{d}{2}}}{c_k^2}\dint_{\R^d}h(y) \varphi(y)\omega_k(y)dy
\nonumber\\ &=& \frac{4^{\gamma + \frac{d}{2}}}{c_k^2}\langle
T_{h\omega_k},\varphi\rangle.
\end{array}
\end{equation}
Thus the relations (\ref{hhhhh}),(\ref{uuuuu}) imply $$S =
T_{h\omega_k}.$$ This relation shows that the support $h$  is
compact. Then $h \in L^2_{k,c}(\R^d)$

%This result is also true on ${\cal S}'(\R^{d})$, more precisely we
%have
 \subsection{ The  Dunkl\
 translation operator and the Dunkl convolution product}

\begin{Def}\hspace*{-2mm}. Let $y \, \in \R^{d}$. The Dunkl translation operator $f \mapsto \tau_y f$ is
defined on ${\cal S}(\R^d)$ by
\begin{equation} \forall \, x \in \R^d, \;
{\cal F}_D (\tau_{y}f)(x)= K(-ix,y){\cal F}_D (f)(y). %
\label{2.38}%
\end{equation}
\end{Def}
\noindent{\bf{Example}} \\ \hspace*{5mm}Let $t > 0$, we have
\begin{equation}\label{tgsz}\forall \, x \in \, \R^d,
\;\tau_{x}(e^{-t||\xi||^2})(y) = \frac{M_k}{t^{\gamma +
\frac{d}{2}}}
 K(\frac{x}{\sqrt{2t}},\frac{y}{\sqrt{2t}})
  e^{-\frac{||x||^2 + ||y||^2}{4t}},
\end{equation}
with $M_k = (2^{\gamma+\frac{d}{2}}c_k)^{-1}$.\\
\noindent{\bf{Remark}}
\\ \hspace*{5mm}  The  operator $\tau_y$, $y
\in \R^d$, can also be  defined    on ${\cal E}(\R^d)$ by
\begin{equation}
\forall \, x \in \R^d, \;\tau_{y}f(x)= (V_k)_x (V_k)_y[(V_k)^{-1}(f)(x+y). %
\label{2.389}%
\end{equation}
(See \cite{T5}).\\ \hspace*{5mm} At the moment an  explicit
formula for the Dunkl
 translation operator is known only in the following two cases. \\
\underline{1$^{st}$ cas }: $d = 1$ and $W = \Z_2$. \\ For all $f
\in C(\R)$ we have $$\begin{array}{ccc}
 \forall \, x \in \R, \tau_{y}f(x) & = &  \frac{1}{2}\dint_{-1}^{1}f(\sqrt{x^2 + y^2 -2xyt})
  (1+\frac{x-y}{\sqrt{x^2 + y^2 -2xyt}})\Phi_k(t)dt\\
   & + &  \frac{1}{2}\dint_{-1}^{1}f(-\sqrt{x^2 + y^2 -2xyt})
  (1-\frac{x-y}{\sqrt{x^2 + y^2 -2xyt}})\Phi_k(t)dt,
\end{array}$$ where $$\Phi_k(t) = \frac{\Gamma(k+\frac{1}{2})}{\sqrt{\pi}%
\Gamma(k)} (1+t)(1-t^2)^{k-1}.$$ Moreover for all $f \in
L^p_k(\R)$, $1 \leq p \leq \infty$, we have
 $$ ||\tau_{y}f||_{k,p} \leq  3 ||f||_{k,p}, \quad 1 \leq
p \leq \infty. $$ (See [10][13]).\\ \underline{2$^{nd}$ cas }: For
all $f \in {\cal E}(\R^d)$ radial we have $$ \forall \, x \in
\R^d, \; \tau_{y}f(x) = V_k [f_0 (\sqrt{||x||^2 + ||y||^2 +2
\langle x,.\rangle })](y),$$ with $f_0$ the function on
$[0,+\infty[$ given by $f(x) = f_0(||x||)$. \\ Moreover for all $f
\in L^p_k(\R^d)$, $1 \leq p \leq \infty$, we have
 $$ ||\tau_{y}f||_{k,p} \leq   ||f||_{k,p}, \quad 1 \leq
p \leq \infty. $$ (See [11][13]).\\ \hspace*{5mm}Using the Dunkl
translation operator, we  define the Dunkl convolution product of
functions as follows (See [11][17]).
\begin{Def}\hspace*{-2mm}. For $f,g$ in $D(\R^d)$, we define the Dunkl convolution
 product  by%
\begin{equation}
\forall \, x \in \R^d, \; f*_{D}g(x)=\int_{\R^d}\tau^{x}f(-y)g(y)d\omega_{k}(y).\label{2.42}%
\end{equation}
\end{Def}
This convolution  is commutative and associative and satisfies the
following properties. (See [13]).\\
\begin{equation}\hspace*{-97mm}i) {\cal F}_D (f*_{D}g) = {\cal
F}_D (f){\cal F}_D (g).\end{equation}
%\begin{Prop}
%{\bf{Young's inequality}}\\
\hspace*{5mm}ii)  Let $1\leq p,q,r\leq+\infty,\;$such that $\frac{1}%
{p}+\frac{1}{q}-\frac{1}{r}=1.\;$If $f\;$is in $L^{p}_{k}(\R^d)$
radial and $g$ an element of $L^{q}_{k}(\R^d),\;$ then $f*_{D}g\;$
belongs to  $L^{r}_{k}(\R^d)\;$ and we have
\begin{equation}
\left\|  f*_{D}g\right\|  _{r,k}\leq \left\|  f\right\| _{p,k}
\left\|  g\right\|  _{q,k}.\label{2.43999}%
\end{equation}
\hspace*{5mm}iii) Let $d = 1$ and $W = \Z_2$. For all
 $f\;$ in $L^{p}_{k}(\R)$ and $g\;$ an element of
$L^{q}_{k}(\R)$, the function $f*_{D}g$ belongs to
$L^{r}_{k}(\R)\;$ with $\frac{1}%
{p}+\frac{1}{q}-\frac{1}{r}=1.\;$ and we have
\begin{equation}
\left\|  f*_{D}g\right\|  _{r,k}\leq 3\left\|  f\right\| _{p,k}
\left\|  g\right\|  _{q,k}.\label{2.43}%
\end{equation}

\section{Functions with compact spectrum } \hspace*{5mm}  First
we recall that the spectrum of a function is the support of its
Dunkl transform. \\We begin this section by the following
definition.
 \begin{Def}\hspace*{-2mm}.
 i) We define the support of $ g \in L^2_k(\R^d)$ and we denote it by
  $\mbox{supp }\, g$,
  the smallest closed  set, outside which the function $g$
 vanishes almost everywhere. \\ \hspace*{5mm}ii) We denote by  $$R_g = \displaystyle \sup_{ \lambda \in suppg}
 ||\lambda||,$$  the radius of the support of $g$.\\
\noindent{{\bf{ Remark}}}\\ \hspace*{5mm} It is clear that $R_g$
is
 finite if and only if, $g$ has compact support.
 \end{Def}
 \begin{Prop}\hspace*{-2mm}. Let $g \in L^2_k(\R^d)$  such that for all $n \in \N$, the
 function
 $||\lambda||^{2n}g(\lambda)$  belongs to $ L^2_k(\R^d)$. Then
\begin{equation}\label{g}
  R_g = \lim_{n \to \infty}\left\{\dint_{\R^d}||\lambda||^{4n}|g(\lambda)|^2
  \omega_k(\lambda)d\lambda\right\}^{\frac{1}{4n}}.
\end{equation}
\end{Prop}
\noindent{\bf{Proof}}\\ \hspace*{5mm} We  suppose  that
$||g||_{k,2} \neq 0$, otherwise $R_g = 0$  and formula (\ref{g})
is trivial.\\ \hspace*{5mm}Assume now that $g$ has compact support
with $R_g
> 0$. Then  $$ \left\{\dint_{\R^d} ||\lambda||^{4n} |g(\lambda)|^2
  \omega_k(\lambda) d\lambda\right\}^{\frac{1}{4n}} \leq
  \left\{\dint_{||\lambda|| \leq R_g}|g(\lambda)|^2
  \omega_k(\lambda)d\lambda\right\}^{\frac{1}{4n}}R_g.$$ Thus we
  deduce that
   $$\limsup_{n \to
\infty}\left\{\dint_{\R^d}||\lambda||^{4n}|g(\lambda)|^2
  \omega_k(\lambda)d\lambda\right\}^{\frac{1}{4n}} \leq
  \limsup_{n \to \infty}\left\{\dint_{||\lambda|| \leq R_g}|g(\lambda)|^2
  \omega_k(\lambda)d\lambda\right\}^{\frac{1}{4n}}R_g = R_g.$$ On the
  other hand, for any positive
  $\varepsilon$ we have
  $$\dint_{R_g - \varepsilon \leq ||\lambda|| \leq R_g}|g(\lambda)|^2
  \omega_k(\lambda)d\lambda > 0.$$ Hence
  $$\liminf_{n \to \infty}\left\{\dint_{\R^d}||\lambda||^{4n}|g(\lambda)|^2
  \omega_k(\lambda)d\lambda\right\}^{\frac{1}{4n}} \geq
 \liminf_{n \to \infty}\left\{\dint_{R_g - \varepsilon \leq ||\lambda|| \leq R_g}
 ||\lambda||^{4n}|g(\lambda)|^2
  \omega_k(\lambda)d\lambda\right\}^{\frac{1}{4n}} \geq R_g -
  \varepsilon.$$
  Thus
$$R_g = \lim_{n \to
\infty}\left\{\dint_{\R^d}||\lambda||^{4n}|g(\lambda)|^2
  \omega_k(\lambda)d\lambda\right\}^{\frac{1}{4n}}.$$
\hspace*{5mm}  We  prove now the assertion in the case where $g$
has
  unbounded support. Indeed For any positive $N$, we have
   $$\dint_{ ||\lambda|| \geq N}|g(\lambda)|^2
  \omega_k(\lambda)d\lambda > 0.$$ Thus
  $$\liminf_{n \to \infty}\left\{\dint_{\R^d}||\lambda||^{4n}|g(\lambda)|^2
  \omega_k(\lambda)d\lambda\right\}^{\frac{1}{4n}} \geq
 \liminf_{n \to \infty}\left\{\dint_{ ||\lambda|| \geq N}
 ||\lambda||^{4n}|g(\lambda)|^2
  \omega_k(\lambda)d\lambda\right\}^{\frac{1}{4n}} \geq N.$$
  This implies that $$\liminf_{n \to \infty}\left\{\dint_{\R^d}||\lambda||^{4n}|g(\lambda)|^2
  \omega_k(\lambda)d\lambda\right\}^{\frac{1}{4n}} = \infty.$$
  \noindent{\bf{Notations.}} We denote by\\
  \hspace*{5mm} - $L^2_{k,R}(\R^d) := \{g \in L^2_{k,c}(\R^d) / R_g =
  R\}$, for $R \geq 0$.\\ \hspace*{5mm} - $D_{R}(\R^d) := \{g \in D(\R^d) / R_g =
  R\}$, for $R \geq 0$.
  \begin{Def}\hspace*{-2mm}.
  We define the Paley-Wiener spaces $PW^2_k(\R^d)$ and  $PW^2_{k,R}(\R^d)$ as
  follows\\
i) $PW^2_k(\R^d)$  is
  the space of  functions $f \in {\cal E}(\R^d)$ satisfying\\
  \hspace*{5mm} a) $\triangle_k^n f \in L^2_{k}(\R^d)$ for all
   $n \in \N$.\\ \hspace*{5mm} b) $R_f^{\triangle_k} := \displaystyle \lim_{n \to  \infty}
  ||\triangle_k ^n f||_{k,2}^{\frac{1}{2n}} < \infty.$\\
ii) $PW^2_{k,R}(\R^d) := \{f \in PW^2_k(\R^d) / R_f^{\triangle_k}
= R\}$.
  \end{Def}

  The real $L^2$-Paley-Wiener theorem for the Dunkl transform can
  be formulated as follows
  \begin{Th}\hspace*{-2mm}.
  The Dunkl transform ${\cal F}_D$ is a bijection \\
  \hspace*{5mm}i) from  $PW^2_{k,R}(\R^d)$ onto
  $L^2_{k,R}(\R^d)$.\\ \hspace*{5mm} ii)from $PW^2_k(\R^d)$
  onto $L^2_{k,c}(\R^d)$,\\
  \end{Th}
\noindent{\bf{Proof}}\\ \hspace*{5mm} i)  Let $g \in
PW^2_{k,R}(\R^d)$. Then from Proposition 2.7 the function ${\cal
F}_D(\triangle_k^n g)(\xi) = (-1)^n ||\xi||^{2n}{\cal F}_D(g)(\xi)
$ belongs to $ L^2_{k}(\R^d)$ for all $n \in \N$. On the other
hand from Theorem 2.3 we deduce that $$ \lim_{n \to
\infty}\left\{\dint_{\R^d}||\xi||^{4n}|{\cal F}_D(g)(\xi)|^2
  \omega_k(\xi)d\xi\right\}^{\frac{1}{4n}} = \lim_{n \to
  \infty}\left\{\dint_{\R^d}|\triangle_k g(x)|^2
  \omega_k(x)dx\right\}^{\frac{1}{4n}} = R.$$
  Thus using Proposition 3.2 we conclude that ${\cal F}_D(g)$ has
  compact support with $ R_{{\cal F}_D(g)} = R$.\\ \hspace*{5mm}
  Conversely let $f \in L^2_{k,R}(\R^d)$. Then $||\xi||^n f(\xi) \in
  L^1_k(\R^d)$ for any  $n \in \N$, and ${\cal F}_D^{-1} f \in
 D(\R^d)$. On the other hand from  Theorem 2.3 we have $$ \lim_{n
\to \infty}\left\{\dint_{\R^d}|\triangle_k^n ({\cal F}_D^{-1}
f)(x)|^2
  \omega_k(x)dx\right\}^{\frac{1}{4n}} = \lim_{n \to
  \infty}\left\{\dint_{\R^d}||\xi||^{4n}|f(\xi)|^2
  \omega_k(\xi)d\xi\right\}^{\frac{1}{4n}} = R.$$
  Thus ${\cal F}_D^{-1}
(f) \in PW^2_{k,R}(\R^d)$.\\  \hspace*{5mm}ii) We deduce ii) from
i).
\begin{Cor}\hspace*{-2mm}.
  The Dunkl transform ${\cal F}_D$ is a bijection
   from  $PW^2_{k}(\R^d)$ onto
  ${\cal H}_{L^2_k}(\C^d)$.
  \end{Cor}
  \noindent{\bf{Proof}}\\ \hspace*{5mm} We deduce the result from
  Theorem 3.4 ii) and Theorem 2.8.
\begin{Def}\hspace*{-2mm}.
 i) The Paley-Wiener space $PW_k(\R^d)$ is
  the space of  functions $f \in {\cal E}(\R^d)$ satisfying\\
  \hspace*{5mm} a) $(1+||x||)^m \triangle_k^n \in L^2_{k}(\R^d)$ for
  all
   $n$,$m$ $\in$ $\N$.\\ \hspace*{5mm} b) $R_f^{\triangle_k} := \lim_{n \to  \infty}
  ||\triangle_k ^n f||_{k,2}^{\frac{1}{2n}} < \infty.$\\
\hspace*{5mm} ii) We have  $PW_{k,R}(\R^d) := \{f \in PW_k(\R^d) /
R_f^{\triangle_k} = R\}$, for $R \geq 0$.
  \end{Def}
  \noindent{{\bf{Remark}}}\\ \hspace*{5mm}
  We notice that the only difference between $PW_k^2(\R^d)$ and
  $PW_k(\R^d)$is the extra requirement of polynomial decay to help
  ensure that ${\cal F}_D (f) \in {\cal E}(\R^d)$.\\

  The real Paley-Wiener theorem for the Dunkl transform of functions in the preceding spaces
  is the following
   \begin{Th}\hspace*{-2mm}.
  The Dunkl transform ${\cal F}_D$ is a bijection \\ \hspace*{5mm}i) from $PW_{k,R}(\R^d)$ onto
  $D_{R}(\R^d)$.\\\hspace*{5mm}ii) from $PW_k(\R^d)$
  onto $D(\R^d)$.\\
  \end{Th}
\noindent{\bf{Proof}}\\ \hspace*{5mm} i)Let $g \in PW_{k,R}(\R^d)
\subset PW^2_{k,R}(\R^d)$. Then ${\cal F}_D (g) \in {\cal
E}(\R^d)$ since $g$ has polynomial decay, and by Theorem 3.4 the
function ${\cal F}_D (g)$ has compact support with $R_{{\cal F}_D
(g)} = R$.\\ \hspace*{5mm} Conversely Let $f \in D_R (\R^d)$, then
${\cal F}_D^{-1}(f) \in {\cal S}(\R^d)$ and ${\cal F}_D^{-1}(f)
\in PW^2_{k,R}(\R^d)$ by Theorem 3.4.\\  \hspace*{5mm}ii) We
deduce the result from the i).
\section{Dunkl transform of functions, with polynomial domain
support} Let $P(x)$ be a non-constant polynomial.
\begin{Th}\hspace*{-2mm}. For any function $f \in {\cal S}(\R^d)$ the following
relation holds
\begin{equation}\label{aze}
 \lim_{n\to\infty} ||P(iT)^n f||_{k,p}^{\frac{1}{n}} = \sup_{y \in supp {\cal F}_D(f)}|P(y)|,
  \;1
  \leq p \leq \infty,
\end{equation}
with $T = (T_1,...,T_d)$.
\end{Th}
\noindent{\bf{Proof}}\\ \hspace*{5mm} We consider $ f \neq 0$ in
${\cal S}(\R^d)$. Set $q = \frac{p}{p-1}$ if $1 < p < \infty$ and
$q = 1$ or $\infty$ if $p = \infty$ or $1$.\\ The proof is divided
in several steps.\\ In the following three steps we suppose that
\begin{equation}\label{wahid}
0 < \sup_{y \in supp \, {\cal
  F}_D(f)}|P(y)| < \infty.\end{equation}
%We see that the
%Dunkl transform of $P(iT)^n f$ is $P^n (\xi) {\cal
%F}_D(f)(\xi).$
\\
{\bf{\underline{First step}}}: In this step we shall prove that
$$\limsup_{n\to\infty} ||P(iT)^n f||_{k,p}^{\frac{1}{n}} \leq
\sup_{y \in supp {\cal F}_D(f)}|P(y)|, \; 1 \leq p \leq \infty.$$
$\bullet$ Let $2
  \leq p < \infty$.
  %We assume that (\ref{l14}) is valid.
  Applying
   Proposition \ref{cher}
   %and using the fact that $${\cal
   %F}_D^{-1}(f)(y) = {\cal F}_D(f)(-y)$$
   we obtain
\begin{eqnarray}\label{l15}
||P(iT)^n f||_{k,p} &\leq& C ||P(\xi)^n {\cal F}_D(f)||_{k,q},\\
\\
& \leq & C (\sup_{y \in supp {\cal F}_D(f)}|P(y)|)^n || {\cal
F}_D(f)
%1_{supp {\cal F}_D(f)}
||_{k,q}. \nonumber
\end{eqnarray}
Thus
\begin{equation}\label{l16}
\limsup_{n\to\infty} ||P(iT)^n f||_{k,p}^{\frac{1}{n}} \leq
\sup_{y \in supp {\cal F}_D(f)}|P(y)|.
\end{equation}
%Let  $\xi_0 \notin U_P$, that means  $P(\xi_0)$. Then there is a
%neighborhood $V_{\xi_0}$ of $\xi_0$ with the property $|P(\xi_0)|
%> \dfrac{1+ |P(\xi_0)| \rangle }{2} > 1$ for $\xi \in
%V_{\xi_0}$. Assume that $p > 1$.  We have
%\begin{eqnarray}\label{l17}\nonumber
%1 &\geq& \limsup_{n\to\infty}||P(\xi)^n {\cal
%F}_D(f)||_{k,q}^{\frac{1}{n}} \geq
%\limsup_{n\to\infty}\{\dint_{V_{\xi_0}}|P(\xi)^n {\cal
%F}_D(f)(\xi)|^q\omega_k(\xi)d\xi\}^{\frac{1}{qn}}\\ \\ &\geq&
%(\dfrac{1+ |P(\xi_0)| }{2})\dint_{V_{\xi_0}}|{\cal
%F}_D(f)(\xi)|^q\omega_k(\xi)d\xi\}^{\frac{1}{qn}}.
%\end{eqnarray}
%Since $(\dfrac{1+ |P(\xi_0)| }{2}) > 1$ and the last limit in
%(\ref{l17}) can be either $1$ or $0$, then
%$$\limsup_{n\to\infty}\{\dint_{V_{\xi_0}}|{\cal
%F}_D(f)(\xi)|^2\omega_k(\xi)d\xi\}^{\frac{1}{qn}} = 0,$$ that
%means $\dint_{V_{\xi_0}}|{\cal F}_D(f)(\xi)|^2\omega_k(\xi)d\xi=
%0$. Thus $\xi_0$ does not belong to the support of ${\cal
%F}_D(f)$. Hence ${\cal F}_D(f) \subset U_P$.\\ Assume now that
%$p=1$. Then
%\begin{eqnarray}\label{l18}\nonumber
%1 &\geq& \limsup_{n\to\infty}||P(\xi)^n {\cal
%F}_D(f)||_{k,\infty}^{\frac{1}{n}} \geq
%\limsup_{n\to\infty}||P(\xi)^n {\cal F}_D(f)||_{L^\infty_k
%(V_{\xi_0}) }^{\frac{1}{n}}\\ \\ &\geq& (\dfrac{1+ |P(\xi_0)|
%}{2})\limsup_{n\to\infty}\limsup_{n\to\infty}|| {\cal
%F}_D(f)||_{L^{\infty}_{k} (V_{\xi_0}) }^{\frac{1}{n}}.
%\end{eqnarray}
%Therefore $|| {\cal F}_D(f)||_{L^\infty_k (V_{\xi_0}) } = 0$, that
%means $\xi_0$ does not belong to the support of ${\cal F}_D(f)$.
%Hence ${\cal F}_D(f) \subset U_P$.\\

%Conversely
\noindent$\bullet$ Suppose now that $1 \leq p < 2$.
% {\cal F}_D(f) \subset U_P$.
  H\"older's inequality gives
\begin{equation}\label{sour}
  ||f||_{k,p}^{p} =
  \dint_{\R^d}(1+||x||^{2})^{-rp}|(1+||x||^{2})^{r}f(x)|^{p}\omega_k(x)dx
  \leq
  ||(1+||x||^{2})^{r}f||_{k,2}^{p}||(1+||x||^{2})^{-rp}||_{k,\frac{2}{2-p}}.
\end{equation}
$$ \leq C  ||(1+||x||^{2})^{r}f||_{k,2}^{p}, $$ for $r > 2\gamma
+d$.\\ Thus, from Proposition 2.7 we obtain $$ ||f||_{k,p}^{p}
\leq C ||(I - \triangle_k )^{r}[{\cal F}_D(f)]||_{k,2}^{p}.$$\\
Consequently for all $n \in \N$, we deduce that
\begin{equation}\label{kopm}
  ||P^n (iT) f||_{k,p} \leq C^{\frac{1}{p}}||(I - \triangle_k
)^{r}[P^n (\xi){\cal F}_D(f)]||_{k,2}.
\end{equation}
On the other hand from Proposition 5.1 of [9] we have, the
following relation:

For all $\mu \in \N^d\backslash\{0\}$ there exist: $t_p ^0, t_p ^1
\in [0,1]$,$p=1,...,|\mu|-1$, such that for all $u \in {\cal
E}(\R^d)$ we have
\begin{eqnarray}\label{65}
 T^{\mu} u(x)  & = & D^{\mu} u(x) + \displaystyle\sum_{\alpha \in R_+
 }\{\displaystyle\sum_{|\beta|=|\mu|}\displaystyle\sum_{p=1}^{|\mu|-1}
Q_{\mu}(t_1 ^0,...,t_p ^0) D^{\beta}u\big(x  - S_{\mu}(t_1 ^0
,...,t_p ^0 )<\alpha,x>\alpha\big) \nonumber
 \\
  & + & \displaystyle\sum_{|\beta'|=|\mu|}
 P_{\mu}(t_1 ^1,...,t_{|\mu|-1}^1)
D^{\beta'}u\big(x  - \widetilde{S}_{\mu}(t_1 ^1
,...,t_{|\mu|-1}^1)<\alpha,x>\alpha\big)
 \},
\end{eqnarray}  where
$Q_{\mu}(t_1,...,t_p),S_{\mu}(t_1,...,t_p)$, $p=1,...,|\mu|$ and
$P_{\mu}(t_1,...,t_{|\mu|-1}),\widetilde{S}_{\mu}(t_1,...,t_{|\mu|-1})$
are polynomials of degree at most $|\mu|$,with respect to each
variable.\\ From this relation and by induction
%and
%the relation (70) in the proof of
%Proposition 5.1
%[9],
one can show that
\begin{equation}\label{htfv}
||(I - \triangle_k )^{r}[P^n (\xi){\cal F}_D(f)(\xi)]||_{k,2} \leq
C n^{2r} ||P^{n-2r}(\xi)\varphi_n(\xi)||_{k,2}, \; n > 2r,
\end{equation}
with $supp \, \varphi_n \subset supp \,{\cal F}_D(f)$ and
$||\varphi_n||_{k,2} \leq C_1 $, where $C_1$ is a constant
independent of $n$.\\ Hence, from the previous inequalities
% and the fact that $|P^{n-2d}(\xi)| \leq 1$
%on the support of $\varphi_n$
we deduce that
\begin{eqnarray}\label{nbvc}\nonumber
 ||P^n (iT) f||_{k,p} &\leq& C^{\frac{1}{p}}n^{2r}
||P^{n-2r}(\xi)\varphi_n(\xi)||_{k,2} \leq C^{\frac{1}{p}}n^{2r}
\sup_{y \in supp {\cal F}_D(f)}|P(y)|^{n-2r}
||\varphi_n(\xi)||_{k,2}\\ \nonumber\\  &\leq& C^{\frac{1}{p}}C_1
n^{2r}\sup_{y \in supp {\cal F}_D(f)}|P(y)|^{n-2r} .
\end{eqnarray}
Thus
\begin{equation}\label{lon}
\limsup_{n\to\infty} ||P(iT)^n f||_{k,p}^{\frac{1}{n}} \leq
\sup_{y \in supp {\cal F}_D(f)}|P(y)|.
\end{equation}
$\bullet$ Let now $p = \infty$. From the relation (22) We have $$
  ||f||_{\infty,k}  \leq   \frac{c_k^2}{4^{\gamma+ \frac{d}{2}}}  ||{\cal
  F}_D(f)||_{k,1}.
  $$
 On the other hand, from Cauchy-Schawrz's inequality we obtain
 $$
 ||{\cal
  F}_D(f)||_{k,1} \leq C_0  ||(1+||\xi||^2)^{\frac{2\gamma+d}{2}}{\cal
  F}_D(f)(\xi)||_{k,2},
 $$
 where $C_0$ is a positive constant.\\
Combining the previous inequalities  and replacing $f$ by $P(iT)^n
f$, we deduce that there exists a positive constant $C$ such that
\begin{equation}\label{hgfd}
 ||P(iT)^n f||_{k,\infty} \leq C ||P^n (\xi)(1+||\xi||^2)^{\frac{2\gamma+d}{2}}{\cal
  F}_D(f)(\xi)||_{k,2}.
\end{equation}
Consequently,
\begin{equation}\label{szaq}
  \limsup_{n \to \infty}  ||P(iT)^n f||_{k,\infty}^{\frac{1}{n}} \leq
  \sup_{y \in supp \, (1+||\xi||^2)^{\frac{2\gamma+d}{2}}{\cal
  F}_D(f)}|P(y)| =  \sup_{y \in supp \, {\cal
  F}_D(f)}|P(y)|.
\end{equation}
Thus from (44), (50) and (52) we have
\begin{equation}\label{pagyz}
  \limsup_{n \to \infty}  ||P(iT)^n f||_{k,p}^{\frac{1}{n}} \leq
   \sup_{y \in supp \, {\cal
  F}_D(f)}|P(y)|, \; 1 \leq p \leq \infty.
\end{equation}
{\bf{ \underline{Second step}}}: In this step  we want to prove
that $$ \lim_{n \to \infty} ||P(iT)^n f||_{k,2}^{\frac{1}{n}} =
 \sup_{y \in supp \, {\cal
  F}_D(f)}|P(y)|.$$
   For any $\varepsilon$, $0 < \varepsilon < \sup_{y \in supp \, {\cal
  F}_D(f)}|P(y)|$, there exists a point $x_0 \in \sup_{y \in supp \, {\cal
  F}_D(f)}|P(y)|$ such that
$$
  |P(x_0)| > \sup_{y \in supp \, {\cal
  F}_D(f)}|P(y)| - \frac{\varepsilon}{2}
$$ As $P$ is a continuous function, there exists a neighborhood
$U_{x_0}$ such that $$
  |P(x)| > \sup_{y \in supp \, {\cal
  F}_D(f)}|P(y)| - \varepsilon, \; x \in U_{x_0}
$$ From Theorem 2.3 we deduce that $$
\begin{array}{lll}
||P(iT)^n f||_{k,2} &=& \frac{c_k^2}{4^{\gamma+
\frac{d}{2}}}||P(\xi)^n {\cal F}_D(f)||_{k,2} \\ &\geq&
\frac{c_k^2}{4^{\gamma+ \frac{d}{2}}} ||P(\xi)^n {\cal
F}_D(f)1_{U_{x_0}}||_{k,2},\end{array} $$ where $1_{U_{x_0}}$ is
the characteristic function of $U_{x_0}$.\\ Thus $$ ||P(iT)^n
f||_{k,2}  \geq \frac{c_k^2}{4^{\gamma+ \frac{d}{2}}}(\sup_{y \in
supp {\cal F}_D(f)}|P(y)| - \varepsilon)^n || {\cal
  F}_D(f) 1_{U_{x_0}}||_{k,2}$$
 This inequality implies,
\begin{equation}\label{usaki}
\liminf_{n \to \infty} ||P(iT)^n f||_{k,2}^{\frac{1}{n}} \geq
(\sup_{y \in supp {\cal F}_D(f)}|P(y)| - \varepsilon) \lim_{n \to
\infty}|| {\cal
  F}_D(f) 1_{U_{x_0}}||_{k,2}^{\frac{1}{n}} =
 \sup_{y \in supp \, {\cal
  F}_D(f)}(|P(y)| - \varepsilon).
\end{equation}
But  $\varepsilon$ can be chosen arbitrarily small, thus from
(\ref{pagyz}) and  (\ref{usaki})   the relation (40) follows for
$p = 2$.\\ \noindent{\bf{ \underline{Third step}}}: In this step
we shall prove that $$ \liminf_{n \to \infty} ||P(iT)^n
f||_{k,p}^{\frac{1}{n}} \geq
 \sup_{y \in supp \, {\cal
  F}_D(f)}|P(y)|,  \;1 \leq p \leq \infty.$$
%Let now $2 < p <
%\infty$. Assume that $supp {\cal F}_D(f) \subset U_P$. Then
%$|P(\xi)| \leq 1$ on the support of ${\cal F}_D(f)$ and therefore,
%by the Proposition \ref{cher} we have
%\begin{equation}\label{nbvca}
 %||P^n (iT) f||_{k,p} \leq
%||P^{n}(\xi){\cal F}_D(f)(\xi)||_{k,q} \leq ||{\cal
%F}_D(f)||_{k,q}.
%\end{equation}
%Consequently, suppose now that (44) is valid.
Since $f \in {\cal S}(\R^d)$, the iteration of the relation (7)
implies the relation
\begin{equation}\label{bhjt}
  \dint_{\R^d} \overline{P^n (-iT) f(x)} P^n (iT) f(x)\omega_k(x)dx
  =  \dint_{\R^d}\overline{f(x)} P^{2n} (iT) f(x)\omega_k(x)dx.
\end{equation}
Hence, by  H\"older's inequality,
\begin{equation}\label{kkkki}
||P^n (iT) f||_{k,2}^2 \leq ||f||_{k,q}||P^{2n} (iT) f||_{k,p}.
\end{equation}
Consequently
\begin{equation}\label{kkkki1}
\lim_{n \to \infty}||P^n (iT) f||_{k,2}^{\frac{1}{n}} \leq
(\lim_{n \to \infty}||f||_{k,q}^{\frac{1}{2n}})\liminf_{n \to
\infty}||P^{2n} (iT) f||_{k,p}^{\frac{1}{2n}} = \liminf_{n \to
\infty}||P^{2n} (iT) f||_{k,p}^{\frac{1}{2n}}.
\end{equation}
Applying now the relation (40)  with $p = 2$, we conclude that
%$supp \, {\cal F}_D(f) \subset U_P$.\\ iii) Let $p = \infty$. Part
%ii) of the proof is still valid with $p = \infty$, $q = 1$.
\begin{equation}\label{ijc}
\sup_{y \in supp \, {\cal
  F}_D(f)}|P(y)| = \lim_{n \to \infty}||P^n (iT)
  f||_{k,2}^{\frac{1}{n}}\leq \liminf_{n \to
\infty}||P^{2n} (iT) f||_{k,p}^{\frac{1}{2n}}.
\end{equation}
 We replace in formula (\ref{kkkki})  the function $f$ by $P(iT) f$ and we
 obtain
\begin{equation}\label{pluyt}
||P^{n+1} (iT) f||_{k,2}^2 \leq ||P(iT)f||_{k,q}||P^{2n+1} (iT)
f||_{k,p}.
\end{equation}
Thus
\begin{equation}\label{tfgr}
  \sup_{y \in supp \, {\cal
  F}_D(f)}|P(y)| = \lim_{n \to \infty}||P^{n+1} (iT)
  f||_{k,2}^{\frac{1}{n+1}}\leq \liminf_{n \to
\infty}||P^{2n+1} (iT) f||_{k,p}^{\frac{1}{2n+1}}.
\end{equation}
Using (\ref{ijc}) and (\ref{tfgr}) we deduce that
\begin{equation}\label{yhwa}
\sup_{y \in supp \, {\cal
  F}_D(f)}|P(y)| \leq \liminf_{n \to
\infty}||P^{n} (iT) f||_{k,p}^{\frac{1}{n}}.
\end{equation}
Then formulas  (\ref{yhwa}) and (53) give  (40). Thus we have
proved the theorem under the condition (\ref{wahid}).\\
\noindent{\bf{ \underline{Fourth step}}}: Suppose now $\sup_{y \in
supp {\cal F}_D(f)}|P(y)| = +\infty.$ Then for any $N
> 0$ there exists a point $x_0 \in supp {\cal F}_D(f)$ such that
$|P(x_0)| \geq 2N$. Since $P$ is a continuous function there
exists a neighborhood $U_{x_0}$ of $x_0$ on which $|P(x)| > N$.
Similarly that the previous calculation of second step  we obtain
$$
\begin{array}{lll} \liminf_{n \to \infty} ||P(iT)^n
f||_{k,2}^{\frac{1}{n}} &\geq& \frac{c_k^2}{4^{\gamma+
\frac{d}{2}}} \liminf_{n \to \infty}||P^n(\xi) {\cal
F}_D(f)1_{U_{x_0}}||_{k,2}^{\frac{1}{n}} ,\\
\\
& \geq &  N \liminf_{n \to \infty}|| f
1_{U_{x_0}}||_{k,2}^{\frac{1}{n}} = N. \end{array}$$ We choose $N$
large, we obtain $$\lim_{n \to \infty} ||P(iT)^n
f||_{k,2}^{\frac{1}{n}} = \infty. $$ Finally if $\sup_{y \in supp
{\cal F}_D(f)}|P(y)| = 0$ the identity (40) is clear   for $p =
2$.\\ Hence the proof of the theorem is finished.
%the inequality (44) follows.\\
\begin{Def}\hspace*{-2mm}.
Let $P$ be a non-constant polynomial and $U_p = \{x \in \R^d, \,
|P(x)| \leq 1 \}$. The set $U_P$ is called a polynomial domain in
$\R^d$.
\end{Def}
\noindent{\bf{Remark}} \\ \hspace*{5mm} A disc is a polynomial
domain. A polynomial domain may be unbounded and nonconvex,  for
example $U = \{x \in \R^d, \, |x_1... x_d| \leq 1 \}$.\\

We have the following result.
\begin{Cor}\hspace*{-2mm}. Let  $f \in {\cal
S}(\R^d)$.
 The Dunkl transform ${\cal F}_D(f)$  vanishes outside a polynomial domain $U_P$, if and only
if,
\begin{equation}\label{l14}
  \limsup_{n\to\infty} ||P(iT)^n f||_{k,p}^{\frac{1}{n}} \leq 1, \; 1
  \leq p \leq \infty.
\end{equation}
\end{Cor}
\noindent{\bf{Remark}} \\ \hspace*{5mm}i) If we take $P(y) = -
||y||^2$, then $P(iT) = \triangle_k$, and Theorem 4.1 and
Corollary 4.3 characterize  functions such that the   support of
their Dunkl transform is a ball.\\  \hspace*{5mm}ii) Theorem 4.1
and Corollary 4.3 generalize also the result obtained in [3].
% The relation (31), but only for functions with
%bounded supports
\section{ Dunkl transform of functions vanishing on a Ball} The
following theorem  gives the radius of the large disc on which the
Dunkl transform of functions in $L^2_{k}(\R^d)$ vanishes  every
where.
\begin{Th}\hspace*{-2mm}.
Let $f \in L^2_{k}(\R^d)$. We consider the sequence
\begin{equation}\label{l}
  f_n(x) = E_n *_D f(x), \; x \in \R^d, \, n \in \N
  \backslash\{0\}.
  %\frac{M_k}{n^{\gamma + \frac{d}{2}}}
  %\dint_{\R^d}K(\frac{x}{\sqrt{2n}},\frac{y}{\sqrt{2n}})
  %e^{-\frac{||x||^2 + ||y||^2}{4n}}f(y)\omega_k(y)dy, \;
  %n=0,1,...,
\end{equation}
where $$ E_n(y) = \frac{c_k}{(4n)^{\gamma +
\frac{d}{2}}}e^{-\frac{||y||^2}{4n}}$$
%with $M_k = (2^{\gamma+\frac{d}{2}}c_k)^{-1}$.\\
Then
\begin{equation}\label{ll}
  \lim_{n \to \infty}\sqrt{-\frac{1}{n}\ln ||f_n||_{k,2}} =
  \lambda_{{\cal F}_D(f)},
\end{equation}
where
\begin{equation}\label{lll}
 \lambda_{{\cal F}_D(f)} = \displaystyle\inf
  \displaystyle\left\{||\xi||, \; \xi \in supp  {\cal
 F}_D(f)\right\}.
\end{equation}\label{TP}
\end{Th}
\noindent{\bf{Remark}}\\ \hspace*{5mm} The function $E_n$ is the
Gauss kernel associated with Dunkl operators. From [11] p. 2424,
we have
\begin{equation}\label{zzzzz}
\forall \, x \in \R^d, \; {\cal F}_D (E_n)(x) = e^{-n||x||^2}.
\end{equation}

\noindent{\bf{Proof of Theorem \ref{TP}}}\\ \hspace*{5mm} First we
remark that from (37) the function $f_n$ is well defined.
%by
%relation (11) and Cauchy-Schwartz inequality  the integral
%$$\dint_{\R^d}K(\frac{x}{\sqrt{2n}},\frac{y}{\sqrt{2n}})
 % e^{-\frac{||x||^2 + ||y||^2}{4n}}f(y)\omega_k(y)dy,$$ is convergent.\\
   We assume that $||f||_{k,2}
> 0
$, otherwise the result is trivial. To prove (\ref{ll}) it is
sufficient to verify the equivalent identity
\begin{equation}\label{llll}
\lim_{n \to \infty} ||f_n||_{k,2}^{\frac{1}{n}} = \exp( -
\lambda_{{\cal F}_D(f)}^2).
\end{equation}
Using (\ref{zzzzz}) and (37) we deduce that the Dunkl transform of
$f_n(x)$ is $ \exp(-n||\xi||^2){\cal F}_D(f)(\xi)$. Then by
applying Theorem 2.3 we obtain
\begin{eqnarray}\label{lllll}
||f_n||_{k,2} &=&
\frac{c_k}{2^{\gamma+\frac{d}{2}}}||\exp(-n||\xi||^2){\cal
F}_D(f)(\xi)||_{k,2} \\ &=&
\frac{c_k}{2^{\gamma+\frac{d}{2}}}||f||_{k,2} \{\dint_{supp {\cal
F}_D(f)}\exp(-2n||\xi||^2)\dfrac{|{\cal
F}_D(f)(\xi)|^2}{||f||_{k,2}^2}
\omega_k(\xi)d\xi\}^{\frac{1}{2}}.\nonumber
\end{eqnarray}
On the other hand it is known that if $m$ is the Lebesque measure
on $\R^d$ and  $U$ a subset of $\R^d$ such that $m(U) = 1$, then
for all $\phi$ in the Lebesgue space $L^p(U,dm)$, $1 \leq p \leq
+\infty$, we have
\begin{equation}\label{l6}
  \lim_{p\to \infty}||\phi||_{L^p(U;dm)} =
  ||\phi||_{L^\infty(U;dm)}.
\end{equation}
By applying  formula (\ref{l6}) with $$U = supp {\cal F}_D(f),
\;\phi = \exp(-||\xi||^2), \; p = 2n, \; \mbox{ and} \;dm(\xi) =
\dfrac{|{\cal F}_D(f)(\xi)|^2}{||f||_{k,2}^2}\omega_k(\xi)d\xi,$$
and using the fact that  $\lim_{n \to +\infty} (\frac{c_k
||f||_{k,2} }{2^{\gamma+\frac{d}{2}}})^{\frac{1}{n}} = 1$.\\ We
obtain
\begin{equation}\label{l7}
\lim_{n \to \infty}||f_n||_{k,2} = \sup_{\xi \in supp  {\cal
F}_D(f)}\exp(-||\xi||^2) = \exp( - \lambda_{{\cal F}_D(f)}^2).
\end{equation}
Which is the relation (\ref{llll}).\\

\hspace*{5mm} A function $f \in L^2_{k}(\R^d)$ is the Dunkl
transform of a function vanishing in a neighborhood of the origin,
if and only if, $\lambda_{{\cal F}_D(f)} > 0$, or equivalently, if
and only if the limit (\ref{llll}) is less than $1$. Thus we have
proved the following result.
\begin{Cor}\hspace*{-2mm}.\label{har}
The condition
\begin{equation}\label{l8}
\lim_{n \to \infty} ||f_n||_{k,2}^{\frac{1}{n}} < 1,
\end{equation}
is necessary and sufficient for a function $f \in L^2_{k}(\R^d)$
to have its Dunkl transform vanishing in a neighborhood of the
origin

\end{Cor}
\noindent{\bf{Remark}}\\ \hspace*{5mm} From Theorem 3.3 and
Corollary \ref{har} it follows that the support of the Dunkl
transform of a function in $L^2_{k}(\R^d)$ is in the tore
$\lambda_{{\cal F}_D(f)} \leq ||\xi|| \leq R_{{\cal F}_D(f)}$, if
and only if,
\begin{equation}\label{l9}
\lambda_{{\cal F}_D(f)} \leq \lim_{n \to
\infty}\sqrt{-\frac{1}{n}\ln ||f_n||_{k,2}}\leq \lim_{n \to
\infty} ||\triangle_k^n f||_{k,2}^{{\frac{1}{2n}}} \leq R_{{\cal
F}_D(f)}.
\end{equation}
\begin{Th}\hspace*{-2mm}. For any function $f \in {\cal S}(\R^d)$ the following
relation holds
\begin{equation}\label{pknfr}
\lim_{n \to \infty} ||\displaystyle\sum_{m = 0}^{\infty}
\frac{(n\triangle_k)^m\, f }{m!}||_{k,p}^{\frac{1}{n}} = \exp( -
\lambda_{{\cal F}_D(f)}^2), \; 1 \leq p \leq \infty.
\end{equation}
In particular, a function $f \in {\cal S}(\R^d)$ is the Dunkl
transform of a function in ${\cal S}(\R^d)$ vanishing in the ball
$B(o,r)$ of center $o$ and radius $r$,  if and only if we have
\begin{equation}\label{pknfr45}
\lim_{n \to \infty} ||\displaystyle\sum_{m = 0}^{\infty}
\frac{(n\triangle_k)^m\, f }{m!}||_{k,p}^{\frac{1}{n}}\leq
\exp(-r^2), \; 1 \leq p \leq \infty.
\end{equation}
\end{Th}
\noindent{\bf{Proof}}\\ \hspace*{5mm} A similar  proof to that of
Theorem 4.1, gives the result.
\section{Dunkl transform of functions, vanishing outside a symmetric body }
\hspace*{5mm} A subset  $K$  of $\R^d$ is called a symmetric body
if $-x \in K$ for all $x \in K$. The set $K^* := \{y \in \R^d, \;
\langle x,y\rangle \leq 1 \; for \, all \, x \in K \}$ is called
the polar set of $K$.
%In \cite{S} a Paley-Wiener theorem a characterize
%functions with spectrum belonging to a symmetric body $K$ has been
%proved.
We state now the following another  real Paley-Wiener theorem.
% for
%functions with symmetric body spectrum.
\begin{Th}\hspace*{-2mm}.
A function $f \in {\cal E}(\R^d)$ is the Dunkl transform of a
function in $L^2_{k}(\R^d)$ vanishing outside a symmetric body
$K$, if and only if, $T^\mu f$ belongs to $L^2_{k}(\R^d)$ for all
$\mu = (\mu_1,...,\mu_d) \in \N^d$, and for all $n \in \N$ we have
\begin{equation}\label{l10}
  \sup_{a \in K^*}||(\langle a,T\rangle)^n f||_{k,2} \leq
  ||f||_{k,2},
\end{equation}
where $T = (T_1,...,T_d)$.
\end{Th}
\noindent{\bf{Proof}}\\ \hspace*{5mm} Let $f \in {\cal E}(\R^d)$
assume $f \neq 0$, otherwise the result is clear. We suppose that
$ {\cal F}_D(f)$ which belongs in $ L^2_{k}(\R^d)$ vanishes out
side a symmetric body $K$. Then $f$ is infinitely differentiable
and belongs to $L^2_{k}(\R^d)$ together with $T^\mu f$ for all
$\mu = (\mu_1,...,\mu_d) \in \N^d$. As the Dunkl transform of $(i
(\langle a,\xi\rangle)^n {\cal F}_D(f)(-\xi)$ is $(\langle
a,T\rangle)^n f$, then by applying Theorem 2.3, we obtain
\begin{equation}\label{l11}
||(\langle a,T\rangle)^n f||_{k,2} =
\frac{c_k}{2^{\gamma+\frac{d}{2}}}||(\langle a,\xi\rangle)^n {\cal
F}_D(f)(\xi)||_{k,2}.
\end{equation}
As $K$  satisfies the symmetric property, we deduce that $|\langle
a,\xi\rangle| \leq 1$ for all $\xi \in K$ and $a \in K^*$. Hence
$$\begin{array}{lll} ||(\langle a,\xi\rangle)^n {\cal
F}_D(f)(.)||_{k,2}^2 &=& \dint_{K}|(\langle a,\xi\rangle)^n {\cal
F}_D(f)(\xi)|^2\omega_k(\xi)d\xi \\ &\leq& \dint_{K}|{\cal
F}_D(f)(\xi)|^2\omega_k(\xi)d\xi =
\frac{4^{\gamma+\frac{d}{2}}}{c_k^2} ||f||_{k,2}^2.\end{array}$$
Thus $$ \sup_{a \in K^*}||(\langle a,T\rangle)^n f||_{k,2} \leq
  ||f||_{k,2}.$$
  \hspace*{5mm} Conversely, we assume that the inequality
  (\ref{l10}) is valid for all $n \in \N$. Since $T^\mu f \in
  L^2_{k}(\R^d)$ for
all $\mu = (\mu_1,...,\mu_d) \in \N^d$. Thus from Proposition 2.7
and Theorem 2.3 and the inequality (\ref{l10}) we obtain for all
$n \in \N$:
\begin{equation}\label{l12}
\sup_{a \in K^*}||(\langle a,\xi\rangle)^n {\cal
F}_D(f)(\xi)||_{k,2} = \frac{2^{\gamma+\frac{d}{2}}}{c_k}\sup_{a
\in K^*}||(\langle a,T\rangle)^n f||_{k,2} \leq
\frac{2^{\gamma+\frac{d}{2}}}{c_k}||f||_{k,2}.
\end{equation}
Let $\xi_0 \notin K$, that means there exists $a_0 \in K^*$ such
that $\langle\xi_0,a\rangle > 1$. Then there is a neighborhood
$U_{\xi_0}$ of $\xi_0$ with the property $\langle\xi,a\rangle >
\dfrac{1+ \langle\xi_0,a\rangle }{2} > 1$, for all $\xi \in
U_{\xi_0}$. Thus for all $n \in \N$:
\begin{eqnarray}\label{l13}\nonumber
\frac{2^{\gamma+\frac{d}{2}}}{c_k}||f||_{k,2} &\geq& \sup_{a \in
K^*}||(\langle a,\xi\rangle)^n {\cal F}_D(f)(\xi)||_{k,2} \geq
(\dint_{U_{\xi_0}}|(\langle a,\xi\rangle)^n {\cal
F}_D(f)(\xi)|^2\omega_k(\xi)d\xi)^{\frac{1}{2}}\\ \\ &\geq&
(\dfrac{1+ \langle\xi_0,a\rangle }{2})^n(\dint_{U_{\xi_0}}|{\cal
F}_D(f)(\xi)|^2\omega_k(\xi)d\xi)^{\frac{1}{2}}. \nonumber
\end{eqnarray}
Since $(\dfrac{1+ \langle\xi_0,a\rangle }{2})^n$ approaches
$\infty$ as $n \to \infty$, (\ref{l13}) holds only if
$$\dint_{U_{\xi_0}}|{\cal F}_D(f)(\xi)|^2\omega_k(\xi)d\xi = 0,$$
this implies that $\xi_0$ does not belongs to the support of
${\cal F}_D(f)$. Hence ${\cal F}_D(f) \subset K$, and Theorem 6.1
is proved.\vspace*{5mm}
\bibliographystyle{unsrt}

\end{document}